\documentclass[12pt,oneside]{article}
\usepackage{amsfonts,graphics}
\usepackage[all,knot,arc,ps,dvips]{xy}

\newcommand{\proof}{{\bf \noindent Proof.\ }}

\def\endpf{\relax\ifmmode\expandafter\endproofmath\else
  \unskip\nobreak\hfil\penalty50\hskip.75em\hbox{}\nobreak\hfil\bull
  {\parfillskip=0pt \finalhyphendemerits=0 \bigbreak}\fi}
\def\bull{\vbox{\hrule\hbox{\vrule\kern3pt\vbox{\kern6pt}\kern3pt\vrule}\hrule}}

\newtheorem{defn}{Definition}[section]
\newtheorem{lemma}[defn]{Lemma}

\newtheorem{theorem}[defn]{Theorem}
\newtheorem{definition}[defn]{Definition}
\newtheorem{remark}[defn]{Remark}
\newtheorem{proposition}[defn]{Proposition}

\newcommand{\zz}{{\mathbb{Z}}}

\newcommand{\nn}{{\mathbb{N}}}
\newcommand{\qq}{{\mathbb{Q}}}

\newcommand{\froyshov}{Fr{\o}yshov }
\newcommand{\ozsvath}{Ozsv\'{a}th }
\newcommand{\szabo}{Szab\'{o} }
\newcommand{\spinc}{\ifmmode{{\rm Spin}^c}\else{${\rm Spin}^c$\ }\fi}
\newcommand{\spinct}{\mathfrak t}
\newcommand{\spincs}{\mathfrak s}
\newcommand{\tors}{{\rm Tors}}

\newcommand{\cals}{\mathcal{S}}
\newcommand{\calt}{\mathcal{T}}
\newcommand{\calh}{\mathcal{H}}
\newcommand{\bj}{\mathbf{j}}

\newenvironment{narrow}[2]{%
 \begin{list}{}{%
  \setlength{\topsep}{0pt}%
  \setlength{\leftmargin}{#1}%
  \setlength{\rightmargin}{#2}%
  \setlength{\listparindent}{\parindent}%
  \setlength{\itemindent}{\parindent}%
  \setlength{\parsep}{\parskip}%
 }%
\item[]}{\end{list}}

\newif\ifpic
\picfalse    
\pictrue   

\begin{document}

\title{Rational homology spheres and four-ball genus}
\author{Brendan Owens and Sa\v{s}o Strle}
\date{}
\maketitle

\begin{abstract}
Using the Heegaard Floer homology of \ozsvath and \szabo we investigate
obstructions to definite intersection pairings bounded by rational homology
spheres.  As an application we obtain new lower bounds for the four-ball genus
of Montesinos links.
\end{abstract}

\section{Introduction}
\label{intro}

Let $Y$ be a rational homology three-sphere and $X$ a smooth negative-definite
four-manifold bounded by  $Y$.  For any $\spinc$ structure $\spinct$ on $Y$ let
$d(Y,\spinct)$ denote the correction term invariant of \ozsvath and \szabo (see
\cite{os4} for the definition; this invariant is the Heegaard Floer homology
analogue of the \froyshov invariant in Seiberg-Witten theory).
It is shown in \cite[Theorem 9.6]{os4} that for each $\spinc$ structure
$\spincs\in\spinc(X)$,
\begin{equation}
\label{thm9.6}
c_1(\spincs)^2+{\rm rk}(H^2(X;\zz))\le4d(Y,\spincs|_Y).
\end{equation}

In order to use this inequality one must study the restriction map
$\spincs\mapsto\spincs|_Y$ from $\spinc(X)$ to $\spinc(Y)$; this
map commutes with the conjugation of \spinc structures. Moreover,
since $\spinc(\,\cdot\,)$ is an affine $H^2(\,\cdot\, ; \zz)$
space, the restriction map is equivariant with respect to the
action of  $H^2(X;\zz)$, where this group acts on $\spinc(Y)$
through the natural group homomorphism $H^2(X;\zz)\to H^2(Y;\zz)$.
In this paper we describe an algorithm that for a given second
Betti number tests each possible four-manifold $X$ (i.e., each
possible intersection form) to see if it can give rise to an
equivariant map for which (\ref{thm9.6}) holds for each
$\spincs\in\spinc(X)$.

The algorithm in principle applies to any rational homology sphere for which
the invariants $d(Y,\spinct)$ are known; this is the case for all Seifert
fibered ones (\cite{os6}; see also \cite{os4} for lens spaces). We describe
the situation in detail for four-manifolds $X$ with $b_2(X) \le 2$. Note that
computations are the simplest for homology lens spaces, since in this case
the number of possible equivariant maps as above is greatly reduced.

We use this algorithm to find obstructions to four-ball genus of a link being
as small as the signature allows it to be. To this end we encode the
information about the link and its slice surface in a manifold pair $(X,Y)$ as
above. Specifically, for a link $L$ in the three-sphere and its slice surface
$F$ in the four-ball, we let $Y$ be the two-fold cover of $S^3$ branched along
$L$, and $X$ be the two-fold cover of $B^4$ branched along $F$; this is
analogous to the slice obstruction of Casson-Gordon \cite{cg} and
Fintushel-Stern \cite{fs}. Applying this to Montesinos links, we get some new
bounds on four-ball genus.

Alternatively, one could try to obtain a lower bound on the
four-ball genus of a knot $K$ by attaching a two-handle to $B^4$
along $K$.  If $K$ is alternating, this approach reproduces the
classical bound given by the signature of $K$; this is reminiscent
of the behaviour of the invariant $\tau(K)$ of \ozsvath and \szabo
\cite{os11}.  This is a purely 3-dimensional invariant defined
using knot Floer homology; it gives the optimal lower bound for
torus knots but agrees with the signature bound for alternating
knots.  By contrast our methods yield new bounds for some
alternating knots.

\

\noindent{\bf Acknowledgements.} Part of this work was completed
while the first author was supported as an EDGE postdoc at
Imperial College, London.  We are grateful to Peter \ozsvath for a
helpful remark.

\section{Four-manifolds bounded by rational homology spheres}
\label{sec:cohomology} In this section we study the relationship
between a  smooth four-manifold  $X$ and its boundary $Y$. We will
assume throughout that $X$ is negative definite. The following is
an extension of \cite[Lemma 3]{cg2}.

\begin{lemma}
\label{lemma:torsion}
Let $Y$ be a rational homology sphere; denote by $h$ the order of $H_1(Y;\zz)$.
Suppose that $Y$ bounds $X$ and denote by $s$ the absolute value of the
determinant of the intersection pairing on $H_2(X,\zz)/\tors$. Then $h=st^2$,
where $st$ is the order of the image of $H^2(X;\zz)$ in $H^2(Y;\zz)$, and $t$
is the order of the image of the torsion subgroup of $H^2(X;\zz)$.
\end{lemma}

\proof
Note that for $b_2(X)>0$, $X$ has a non-degenerate integer intersection form
$$Q_X\colon H_2(X;\zz)/\tors \otimes H_2(X;\zz)/\tors \longrightarrow \zz;$$
we denote the absolute value of the determinant of this pairing by $s$.  If
$b_2(X)=0$, then set $s=1$. The long exact sequence of the pair $(X,Y)$ yields
the following (with integer coefficients):

\addtolength{\arraycolsep}{-3pt}
$$\begin{array}{cccccccccccc}
0\longrightarrow&H^2(X,Y)&\stackrel{j^*}{\longrightarrow}&H^2(X)&\longrightarrow&H^2(Y)
&\longrightarrow&H^3(X,Y)&\longrightarrow&H^3(X)&\longrightarrow0,\\ &\|&&\|&&
&&\|&&\|&\\ &\zz^b\oplus T_2&&\zz^b\oplus T_1&& && T_1 && T_2 &
\end{array}$$
\addtolength{\arraycolsep}{+3pt}

\noindent where $T_1, T_2$ are torsion groups, and $b=b_2(X)$ (we
may assume that $b_1(X)=0$; if not one may surger out $b_1$
without changing the conclusion of the lemma). With respect to
appropriate bases for (a compatible choice of) free parts of
$H^2(X,Y)$ and $H^2(X)$, we have
$$j^*=\left[\matrix{Q & 0 \cr * & \tau}\right],$$
where $Q$ is the matrix representation of the intersection pairing on
$H_2(X;\zz)/\tors$. Note that $\tau\colon T_2 \to T_1$ is a monomorphism;  let
$t=|T_1|/|T_2|$. It
follows that $h=st^2$, as $Q$ can be thought of as a presentation matrix for a
group of order $s$.
\endpf

To state the basic relation between $X$ and $Y$ more explicitly, we need to
understand the restriction map from \spinc structures on $X$ to those on $Y$.
Let $\calt$ be the image of the torsion subgroup of $H^2(X;\zz)$ in
$\calh:=H^2(Y;\zz)$, and let $\cals$ be the quotient of $H^2(X;\zz)$ by
the sum of its torsion subgroup and the image of $H^2(X,Y;\zz)$. After fixing
affine identifications of $\spinc(\,\cdot\,)$ with $H^2(\,\cdot\,;\zz)$, the
restriction map from $\spinc(X)$ to $\spinc(Y)$ induces an affine monomorphism
$$\rho : \cals \to \calh/\calt.$$
For appropriate choices of origins in the spaces of \spinc
structures, $\rho$ becomes a group homomorphism and the
conjugation of \spinc structures, denoted by $\bj$, corresponds to
multiplication by $-1$. Choose an identification $\spinc(Y) \cong
\calh$ so that a spin structure corresponds to $0 \in \calh$, and
let $0 \in \cals$ correspond to the class of a \spinc structure on
$X$ whose Chern class belongs to the sum of the torsion subgroup
of $H^2(X;\zz)$ and the image of $H^2(X,Y;\zz)$. If the order of
$\calh$ is odd then there is a unique $\bj$-fixed element in each
of $\cals$ and $\calh/\calt$ and $\rho$ is a group homomorphism.
In general, any $\bj$-fixed element (i.e., any element of order 2)
can be used as origin; to make $\rho$ a group homomorphism one
needs to choose the right spin structure on $Y$.

We define two (rational-valued) functions on $\cals$; one induced by the intersection
pairing on $X$ and the other coming from the correction term on $Y$.
For each $\alpha \in \cals$ let $sq(\alpha)$ be the largest square of the Chern
class of any \spinc structure on $X$ in the equivalence class $\alpha$, and let
$d_\rho(\alpha)$ be the minimal value of the correction term for $Y$ on the
coset $\rho(\alpha)$.

\begin{theorem}
\label{thm:rhomap}
Suppose that a rational homology
sphere $Y$ bounds a negative definite manifold $X$. Then, with above notation,
$$ sq(\alpha) +b_2(X) \le 4d_\rho(\alpha)$$
for all $\alpha \in \cals$.
\end{theorem}
\begin{remark}
Since both sides in the above inequality are $\bj$-invariant, one may work over $\cals/\bj$.
\end{remark}

\proof
This follows from \cite[Theorem 9.6]{os4} and the fact that changing a \spinc
structure on $X$ by a torsion line bundle does not change its square.
\endpf

We note that if $Y$ is a homology lens space, we can choose a labelling
$\{\spinct_j : j=0,\ldots, h-1\}$ of \spinc structures on $Y$ corresponding to an
isomorphism $\calh\cong\zz/h$. Similarly, we label a set of \spinc
structures $\{\spincs_i : i=0,\ldots, s-1\}$ on $X$, where $\spincs_i$ has maximal
square in its equivalence class $i \in \zz/r\cong \cals$; here $s$ denotes
the absolute value of the determinant of the intersection pairing on
$H_2(X;\zz)$. We call such a collection of \spinc structures on $X$ an {\em
optimal set of \spinc structures}. The condition of Theorem \ref{thm:rhomap}
can then be expressed as follows:  for any $i=0,\ldots,s-1$
$$c_1(\spincs_i)^2+b_2(X)\le 4d(Y,\spinct_{\rho(i)+kst}) \hbox{\ \ for all\ \ }
k=0,\ldots,t-1.$$

\section{Application to links}
\label{sec:knots}
Let $L$ be an oriented link with $\mu$ components in the three-sphere; denote
its signature by $\sigma(L)$.  The unlinking number (or unknotting number)
$u(L)$ is
the minimal number of crossing changes in any diagram of $L$ which yield the
trivial $\mu$-component link.

The four-ball genus $g^*(L)$ of $L$ is defined
to be the minimal genus of a (connected) oriented surface $F$ admitting a
smooth embedding into $B^4$ which maps $\partial F$ to $L$.  An easy argument
shows that $g^*(L)\le u(L)$.  A classical result
due to Murasugi \cite{m} states that
\begin{equation}
\label{eqn:murasugi}
g^*(L)\ge \frac{|\sigma(L)|-\mu+1}{2}.
\end{equation}
Suppose that this bound is attained and fix such a connected minimal surface
$F$. Let $X$ be the branched double cover of $B^4$ along $F$. Then
$b_1(X)=0,\ b_2(X)=2g^*(L)+\mu-1$, and the signature of $X$ is given by
$\sigma(L)$ (\cite{kt}). After possibly changing its orientation, we may assume
that $X$ is negative-definite.

Note that $Y=\partial X$  is the double cover of $S^3$ branched along $L$. If
$Y$ is a rational homology sphere (which is the case if the determinant
$h=|\Delta_L(-1)|$ of $L$ is non-zero; in this case $h$ is the order of
$H_1(Y;\zz)$), we may apply Theorem \ref{thm:rhomap}. We will spell this out in
more detail in Section \ref{sec:algorithm}.

In Section \ref{sec:examples} we list some resulting bounds on the four-ball
genus of Montesinos links. In the rest of this section we discuss other
classical bounds on the four-ball genus; we describe Montesinos links, their
double branched covers, and a spanning surface;
and we recall the formulas from \cite{os4,os6} for the
correction term of Seifert fibered rational homology spheres.

\subsection{Bounds on four-ball genus from Seifert matrices}
\label{subsec:taylor}
In the case of a knot $K$, the signature in
(\ref{eqn:murasugi}) may be replaced by any Tristram-Levine
signature $\sigma_\omega(K)$, where $\omega \in S^1-1$, yielding
potentially stronger bounds.  These signatures may be computed
from any Seifert matrix associated to $K$.  A stronger bound is
given by Taylor \cite{t}, which we now describe.

Let $M \in \zz^{a\times a}$ be any Seifert matrix for $K$. Then
$M$ defines a pairing $\lambda$ on $\zz^a$ by $\lambda(x,y)=x M
y^T$. Denote by $z(M)$ the maximal rank of a nullspace of
$\lambda$, that is a sublattice $N$ such that $\lambda(x,y)=0$ for
all $x,y\in N$.  Taylor defines an invariant $m(K)=a/2-z(M)$, and
he proves the following inequalities for any $\omega$:
\begin{equation}
\label{eqn:taylor}
g^*(K)\ge m(K) \ge \frac{|\sigma_\omega(K)|}{2}.
\end{equation}
In Section \ref{sec:examples} we will provide examples of knots $K$ with
$m(K)=1$ but for which it follows from Theorem \ref{thm:rhomap} that
$g^*(K)>1$.

\subsection{Montesinos links and Seifert fibered spaces}
\label{subsec:mont}

For more details on Montesinos links and their classification see \cite{bz}. In
Definitions \ref{def:mont} and \ref{def:seifert}, $e$ is any integer and
$(\alpha_1,\beta_1),(\alpha_2,\beta_2),\ldots,(\alpha_r,\beta_r)$ are coprime
pairs of integers, with $\alpha_i>1$.

\begin{definition}
\label{def:mont}
A Montesinos link
$M(e;(\alpha_1,\beta_1),(\alpha_2,\beta_2),\ldots,(\alpha_r,\beta_r))$ is a
link which has a projection as shown in Figure \ref{fig:mont}(a).  There are
$e$ half-twists on the left side.  A box \framebox{$\alpha,\beta$} represents a
{\em rational tangle of slope} $\alpha/\beta$:  given a continued fraction
expansion
$$\frac{\alpha}{\beta}=[a_1,a_2,\ldots,a_m]:=a_1-\frac{1}{a_2-\raisebox{-3mm}{$\ddots$
\raisebox{-2mm}{${-\frac{1}{\displaystyle{a_m}}}$}}}\, ,$$
the rational tangle of slope $\alpha/\beta$ consists of the four string braid
$\sigma_2^{a_1}\sigma_1^{a_2}\sigma_2^{a_3}\sigma_1^{a_4}\ldots\sigma_i^{a_m}$,
which is then closed on the right as in Figure \ref{fig:mont}(b) if $m$ is odd
or (c) if $m$ is even.
\end{definition}

\begin{figure}[htbp]
\begin{center}
\ifpic
\leavevmode
\xygraph{
!{0;/r1.0pc/:}
!{*{(a)}}[rd]
!{\vcap[4]}[r]
!{\vcap[2]}[l]
!{\xcapv@(0)}[ur]
!{\xcapv@(0)}[l]
!{\xcapv@(0)}[ur]
!{\xcapv@(0)}[l]
!{\vtwist}
!{\vtwist}
!{\vtwist}
!{\xcapv@(0)}[ur]
!{\xcapv@(0)}[l]
!{\xcapv@(0)}[ur]
!{\xcapv@(0)}[l]
!{\vcap[-4]}[r]
!{\vcap[-2]}[r(2.5)u(6.3)]
!{*+[F]{\alpha_1,\beta_1}}[l(0.5)d(0.75)]
!{\xcapv[0.5]@(0)}[ur]
!{\xcapv[0.5]@(0)}[l(0.5)d(0.2)]
!{*+[F]{\alpha_2,\beta_2}}[l(0.5)d(0.75)]
!{\xcapv[0.5]@(0)}[ur]
!{\xcapv[0.5]@(0)}[l(0.5)u(0.1)]
!{*{.}}[d(0.2)]
!{*{.}}[d(0.2)]
!{*{.}}[d(0.35)l(0.5)]
!{\xcapv[0.5]@(0)}[ur]
!{\xcapv[0.5]@(0)}[l(0.5)d(0.2)]
!{*+[F]{\alpha_r,\beta_r}}[r(5)u(7.25)]
!{*{(b)}}[r]
!{\xbendd[-2]@(0)}[u]
!{\xbendd-@(0)}[dl]
!{\xcaph[1]@(0)}[l(2)d]
!{\xbendu[-2]@(0)}[d]
!{\xbendu-@(0)}[ul]
!{\xcaph[1]@(0)}[l(0.5)u(2)]
!{\xcaph[1]@(0)}[ld]
!{\xcaph[1]@(0)}[ld]
!{\xcaph[1]@(0)}[ld]
!{\xcaph[1]@(0)}[uuu]
!{\xcaph[1]@(0)}[ld]
!{\htwist}[ldd]
!{\xcaph[1]@(0)}[uuu]
!{\xcaph[1]@(0)}[ld]
!{\htwist|{a_1}}[ldd]
!{\xcaph[1]@(0)}[uuu]
!{\xcaph[1]@(0)}[ld]
!{\htwist}[ldd]
!{\xcaph[1]@(0)}[uuu]
!{\htwistneg>{\,\,\,\,\,\,\,\,\,\,\,\,\,\,\,\,\,\,\,\,a_2}}[ldd]
!{\xcaph[1]@(0)}[ld]
!{\xcaph[1]@(0)}[uuu]
!{\htwistneg}[ldd]
!{\xcaph[1]@(0)}[ld]
!{\xcaph[1]@(0)}[uuu]
!{\xcaph[1]@(0)}[ld]
!{\htwist|{a_3}}[ldd]
!{\xcaph[1]@(0)}[uuu]
!{\hcap}[dd]
!{\hcap}[uu]
[r(3)u]
!{*{(c)}}[r]
!{\xbendd[-2]@(0)}[u]
!{\xbendd-@(0)}[dl]
!{\xcaph[1]@(0)}[l(2)d]
!{\xbendu[-2]@(0)}[d]
!{\xbendu-@(0)}[ul]
!{\xcaph[1]@(0)}[l(0.5)u(2)]
!{\xcaph[1]@(0)}[ld]
!{\xcaph[1]@(0)}[ld]
!{\xcaph[1]@(0)}[ld]
!{\xcaph[1]@(0)}[uuu]
!{\xcaph[1]@(0)}[ld]
!{\htwist}[ldd]
!{\xcaph[1]@(0)}[uuu]
!{\xcaph[1]@(0)}[ld]
!{\htwist|{a_1}}[ldd]
!{\xcaph[1]@(0)}[uuu]
!{\xcaph[1]@(0)}[ld]
!{\htwist}[ldd]
!{\xcaph[1]@(0)}[uuu]
!{\htwistneg}[ldd]
!{\xcaph[1]@(0)}[ld]
!{\xcaph[1]@(0)}[uuu]
!{\htwistneg>{a_2}}[ldd]
!{\xcaph[1]@(0)}[ld]
!{\xcaph[1]@(0)}[uuu]
!{\htwistneg}[ldd]
!{\xcaph[1]@(0)}[ld]
!{\xcaph[1]@(0)}[uuu]
!{\hcap[3]}[d]
!{\hcap}
}
\else \vskip 5cm \fi
\begin{narrow}{0.4in}{0.2in}
\caption{
{\bf{Montesinos links and rational tangles.}}
Note that $e=3$ in (a).  Also (b) and (c) are both
representations of the rational tangle of slope 10/3:
$$10/3=[3,-2,1]=[3,-3]$$ (and one can switch between (b) and (c) by simply
moving the last crossing).}
\label{fig:mont}
\end{narrow}
\end{center}
\end{figure}

A two-bridge link $S(p,q)$ (or rational link, or 4-plat) is the reflection of
the link formed by closing the
rational tangle \framebox{$p,q$} with two trivial bridges.  This is
equal to the Montesinos link $M(e;(q,eq+p))$ for any $e$.

\begin{figure}[htbp]
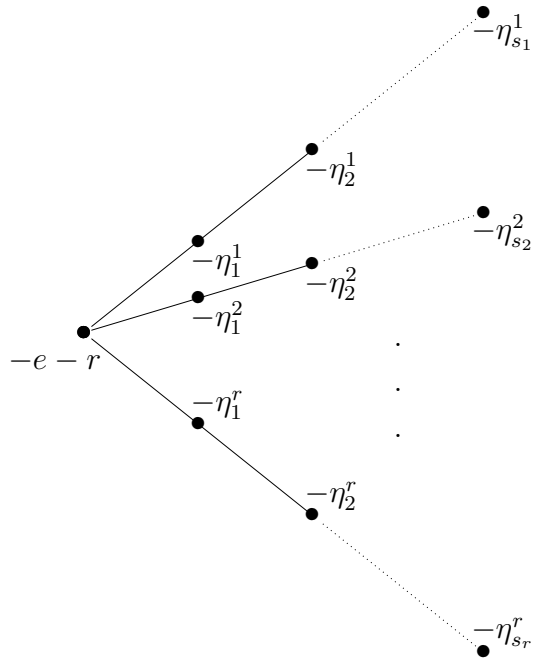

  \begin{center}
\ifpic
    \leavevmode
    \[ \xy/r1.8pc/:
  (0,0)="centre", "centre"+/dl3ex/*{-e-r},
  (2,1.6)="1,1", "1,1"+/dr2ex/*{-\eta_1^1},
  (4,3.2)="1,2", "1,2"+/dr2ex/*{-\eta_2^1},
  (7,5.6)="1,s", "1,s"+/dr2ex/*{-\eta_{s_1}^1},
  (2,0.6)="2,1", "2,1"+/dr2ex/*{-\eta_1^2},
  (4,1.2)="2,2", "2,2"+/dr2ex/*{-\eta_2^2},
  (7,2.1)="2,s", "2,s"+/dr2ex/*{-\eta_{s_2}^2},
  (2,-1.6)="r,1", "r,1"+/ur2ex/*{-\eta_1^r},
  (4,-3.2)="r,2", "r,2"+/ur2ex/*{-\eta_2^r},
  (7,-5.6)="r,s", "r,s"+/ur2ex/*{-\eta_{s_r}^r},
  "centre"*{\bullet};      
  "1,2" **\dir{-},         
  "1,1"*{\bullet};         
  "1,2"*{\bullet};         
  "1,s" **\dir{.},         
  "1,s"*{\bullet};         
  "centre"*{\bullet};      
  "2,2" **\dir{-},         
  "2,1"*{\bullet};         
  "2,2"*{\bullet};         
  "2,s" **\dir{.},         
  "2,s"*{\bullet};         
  "centre"*{\bullet};      
  "r,2" **\dir{-},         
  "r,1"*{\bullet};         
  "r,2"*{\bullet};         
  "r,s" **\dir{.},         
  "r,s"*{\bullet};         
  (5.5,-0.2)*{.};      
  (5.5,-1)*{.};
  (5.5,-1.8)*{.}
\endxy \]
\else \vskip 5cm \fi
    \caption{{\bf Plumbing description of Seifert fibered space.}}
    \label{fig:graph}
  \end{center}
\end{figure}

\begin{definition}
\label{def:seifert}
The Seifert fibered space
$Y(e;(\alpha_1,\beta_1),(\alpha_2,\beta_2),\ldots,(\alpha_r,\beta_r))$
is the oriented boundary of the four-manifold obtained by plumbing
disk bundles over
the two-sphere according to the weighted graph shown in Figure
\ref{fig:graph}.  To each vertex $v$ with multiplicity $m(v)$, associate a
disk bundle over $S^2$ with Euler number $m(v)$.  The bundles associated to
two vertices are plumbed precisely when the vertices are connected by an
edge.
(See \cite{hnk,o} for details on plumbing.)
The multiplicities on the graph are obtained from continued fraction
expansions
$$\frac{\alpha_i}{\alpha_i-\beta_i}=[\eta_1^i,\eta_2^i,\ldots,\eta_{s_i}^i].$$
\end{definition}

A lens space $L(p,q)$ is a special case of the above; it is the boundary of
the plumbed four-manifold associated to a linear graph with weights
$-a_1,-a_2,\ldots,-a_m$, where
${\displaystyle  \frac{p}{q}}=[a_1,a_2,\ldots,a_m]$.
This is equal
to the Seifert fibered space $Y(-e;(q,eq+p))$ for any $e$.

A Seifert fibered space
$Y(e;(\alpha_1,\beta_1),\ldots,(\alpha_r,\beta_r))$
is a rational homology sphere if and only if
$\displaystyle{e+\sum_{i=1}^r\frac{\beta_i}{\alpha_i}\ne 0}$.

\begin{proposition}
\label{prop:2cover}
The branched double cover of $S^3$ along the Montesinos link\linebreak
$M(e;(\alpha_1,\beta_1),\ldots,(\alpha_r,\beta_r))$ is
the Seifert fibered space 
$Y(-e;(\alpha_1,\beta_1),\ldots,(\alpha_r,\beta_r))$.
\end{proposition}

Note it follows from Proposition \ref{prop:2cover} that the branched
double cover of $S(p,q)$ is $L(p,q)$.
\newline

\proof
The original proof is in \cite{mont}.  The result is also proved in \cite{bz}
but note that on p.~197 an $e$-twist should correspond to
$\alpha_0=1,\beta_0=-e$ (rather than $\beta_0=e$).  Since it is particularly
important that we correctly identify the branched cover as an oriented
manifold we will sketch a proof here.

\begin{figure}[htbp]
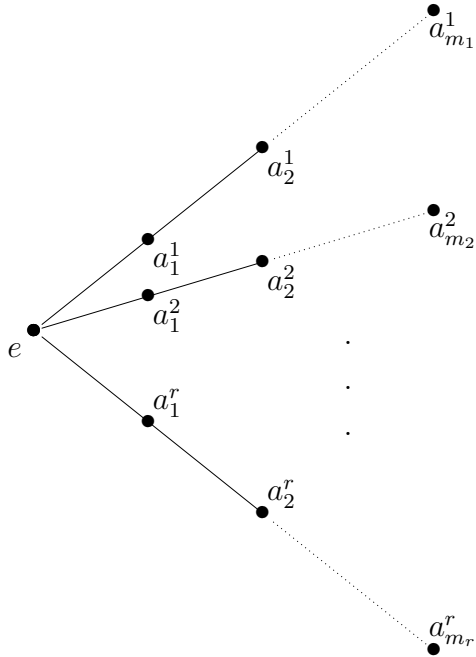

  \begin{center}
\ifpic
    \leavevmode
    \[ \xy/r1.8pc/:
  (0,0)="centre", "centre"+/dl2ex/*{e},
  (2,1.6)="1,1", "1,1"+/dr2ex/*{a_1^1},
  (4,3.2)="1,2", "1,2"+/dr2ex/*{a_2^1},
  (7,5.6)="1,s", "1,s"+/dr2ex/*{a_{m_1}^1},
  (2,0.6)="2,1", "2,1"+/dr2ex/*{a_1^2},
  (4,1.2)="2,2", "2,2"+/dr2ex/*{a_2^2},
  (7,2.1)="2,s", "2,s"+/dr2ex/*{a_{m_2}^2},
  (2,-1.6)="r,1", "r,1"+/ur2ex/*{a_1^r},
  (4,-3.2)="r,2", "r,2"+/ur2ex/*{a_2^r},
  (7,-5.6)="r,s", "r,s"+/ur2ex/*{a_{m_r}^r},
  "centre"*{\bullet};      
  "1,2" **\dir{-},         
  "1,1"*{\bullet};         
  "1,2"*{\bullet};         
  "1,s" **\dir{.},         
  "1,s"*{\bullet};         
  "centre"*{\bullet};      
  "2,2" **\dir{-},         
  "2,1"*{\bullet};         
  "2,2"*{\bullet};         
  "2,s" **\dir{.},         
  "2,s"*{\bullet};         
  "centre"*{\bullet};      
  "r,2" **\dir{-},         
  "r,1"*{\bullet};         
  "r,2"*{\bullet};         
  "r,s" **\dir{.},         
  "r,s"*{\bullet};         
  (5.5,-0.2)*{.};      
  (5.5,-1)*{.};
  (5.5,-1.8)*{.}
\endxy \]
\else \vskip 5cm \fi
    \caption{{\bf Plumbing description of Montesinos link.}}
    \label{fig:montgraph}
  \end{center}
\end{figure}

We start with an alternative description of the Montesinos link\linebreak
$M(e;(\alpha_1,\beta_1),\ldots,(\alpha_r,\beta_r))$.  For
each $i$, let
$$\frac{\alpha_i}{\beta_i}=[a_1^i,a_2^i,\ldots,a_{m_i}^i].$$
Then $M(e;(\alpha_1,\beta_1),\ldots,(\alpha_r,\beta_r))$
is obtained by plumbing twisted bands according to the graph in
Figure \ref{fig:montgraph}, and then taking the boundary.
Here each vertex represents a twisted band, that
is a $D^1$-bundle over $S^1$, embedded in $S^3$, with the number of half-twists
given by the multiplicity of the vertex.  For example, Figure
\ref{fig:mont}(a) is (the boundary of) a band with 3 half-twists, if $r=0$.
Bands are plumbed
together precisely when the corresponding vertices are adjacent.

We now want to describe the double branched cover of such a plumbed link.
Start with the case of a single vertex, with weight $a$.
This gives the two-bridge link $L=S(a,1)$ formed by
closing the four-string braid $\sigma_2^a$.
Split $S^3$ along a 2-sphere which separates the link into 4 arcs, so that
the braid is contained in one component of $S^3-S^2$. (If $L$ is pictured
as in Figure \ref{fig:mont}(b), the 2-sphere may be drawn as a vertical line
through $L$ on one side of the twists.)
This gives $(S^3,L)$ as a union of two balls, each containing two arcs.
The branched double covers of these are solid tori, which inherit
an orientation from $S^3$.  Choose a meridian and longitude pair with
intersection number $+1$ on the boundary of each torus.
With respect to these ordered bases, the map induced on homology by the gluing
map on the boundary of the solid tori is represented by a $2\times2$ matrix.

The braid operation $\sigma_2$ lifts to a right-handed Dehn twist about the
longitude of either torus (choose one).  This has matrix
$\left(\begin{array}{cc}
1 & 0\\
-1 & 1 \end{array} \right)$
in the chosen basis for that torus.
Composing the Dehn twists and  changing
basis to that of the other torus yields the matrix product
$$\left(\begin{array}{cc}
-1 & 0\\
0 & 1 \end{array} \right)
\left(\begin{array}{cc}
1 & 0\\
-a & 1 \end{array} \right)=
\left(\begin{array}{cc}
-1 & 0\\
-a & 1 \end{array} \right).$$
This is precisely the gluing map for the circle bundle over $S^2$ with Euler
number $a$.

Now consider a graph with two vertices labeled $a_1,a_2$ which are
joined by an edge.  The resulting plumbed link is equivalent to
the two-bridge link $L$ formed by closing the four string braid
$\sigma_2^{a_1}\sigma_1\sigma_2\sigma_1\sigma_2^{a_2}$.  As above
split $S^3$ along a 2-sphere to one side of the braid.  The braid
$\sigma_1$ lifts to a right-handed Dehn twist about the meridian,
with matrix $\left(\begin{array}{cc} 1 & 1\\ 0 & 1 \end{array}
\right)$.  Thus the double branched cover of $L$ is the union of
two solid tori with the gluing map given by the product
\begin{eqnarray*}
&&
\left(\begin{array}{cc}
-1 & 0\\
0 & 1 \end{array} \right)
\left(\begin{array}{cc}
1 & 0\\
-a_1 & 1 \end{array} \right)
\left(\begin{array}{cc}
1 & 1\\
0 & 1 \end{array} \right)
\left(\begin{array}{cc}
1 & 0\\
-1 & 1 \end{array} \right)
\left(\begin{array}{cc}
1 & 1\\
0 & 1 \end{array} \right)
\left(\begin{array}{cc}
1 & 0\\
-a_2 & 1 \end{array} \right)\\
&=&
\left(\begin{array}{cc}
-1 & 0\\
-a_1 & 1 \end{array} \right)
\left(\begin{array}{cc}
0 & 1\\
1 & 0 \end{array} \right)
\left(\begin{array}{cc}
-1 & 0\\
-a_2 & 1 \end{array} \right),
\end{eqnarray*}
which is the gluing map for the boundary of the manifold formed by plumbing
together disk bundles over $S^2$ with Euler numbers $a_1,a_2$.

It is now not hard to see that in general if $L$ is the plumbed link
associated to a weighted tree $T$ then the double branched cover of $L$ is the
Seifert fibered space associated to $T$.  According to
Definition \ref{def:seifert}, the Seifert fibered space obtained from
the graph in Figure \ref{fig:montgraph} is
$Y(-e-r;(\alpha_1,\alpha_1+\beta_1),\ldots,(\alpha_r,\alpha_r+\beta_r))$.
Note that this is orientation-preserving diffeomorphic to
$Y(-e;(\alpha_1,\beta_1),\ldots,(\alpha_r,\beta_r))$ as claimed.\endpf

\begin{remark}
\label{rmk:orient}
We have used the same orientation convention as Orlik \cite{o} and\linebreak
Hirzebruch-Neumann-Koh \cite{hnk} for lens spaces and Seifert fibered
spaces.  However the opposite convention for lens spaces is used in \cite{os4}.
\end{remark}

\subsection{A spanning surface for Montesinos links}
\label{subsec:spanning}
We describe an orientable spanning surface $\Sigma$ in $S^3$ for the
Montesinos link $L=M(e;(\alpha_1,\beta_1),\ldots,(\alpha_r,\beta_r))$ which is
a generalisation of that shown in \cite[12.26]{bz} for 2-bridge links (see also
\cite{b}). For knots this will enable us to compute the signature, and also in
some cases the Taylor invariant $m(K)$.  For links with more than one component
both the signature and the four-ball genus depend on a choice of orientation;
we will choose the orientation given by $L=\partial\Sigma$ (for either
orientation of $\Sigma$).

Note the following equivalence of unoriented links
\begin{equation}
\label{eqn:mont} M(e;\ldots,(\alpha_i,\beta_i),\ldots) =
M(e+1;\ldots,(\alpha_i,\alpha_i+\beta_i),\ldots),
\end{equation}
which follows from Proposition \ref{prop:2cover} and the fact that
Montesinos links are classified, up to orientation, by their
branched double covers.  Fixing the surface $\Sigma$ will require
fixing a choice of invariants
$(e;(\alpha_1,\beta_1),\ldots,(\alpha_r,\beta_r))$.

Recall from \cite[12.16]{bz} that if $\alpha,\beta$ are coprime
with $\beta$ odd then we may choose a continued fraction expansion
$$\frac{\alpha}{\beta}=[a_1,a_2,\ldots,a_m]$$
with $m$ odd and $a_2,a_4,\ldots$ even.  Similarly if $\alpha$ is
odd one may choose an expansion with $m,a_1,a_3,\ldots$ even.

Colour black or white, in chessboard fashion, the regions of $S^2$ that form
the complement of the projection in Figure \ref{fig:mont}(a).  Start by
colouring black the twisted band on the left.  There are then two cases to
consider.

\noindent {\bf Case 1:} $\alpha_i$ is odd for $i=1,\ldots,r$.
Assume, using (\ref{eqn:mont}) if necessary, that
$$1\le\beta_i<\alpha_i \hbox{\ \ for all\ \ } i=1,\ldots,r.$$
Then for each $i$, choose a continued fraction expansion
$$\frac{\alpha_i}{\beta_i}=[a_1^i,a_2^i,\ldots,a_{m_i}^i]$$
with $m_i, a_1^i,a_3^i,\ldots,a_{m_i-1}^i$ even.  The white surface is
orientable in the resulting diagram.

\noindent {\bf Case 2:} $\{\alpha_i\}$ are not all odd.  Using
(\ref{eqn:mont}) we may assume each $\beta_i$ is the smallest
positive odd integer in its congruence class mod $\alpha_i$. We
also require that $e\equiv r\pmod2$.  If this does not hold,
choose $(\alpha_j,\beta_j)$ such that $\alpha_j$ is even and
$\displaystyle{\frac{\beta_j}{\alpha_j}=
\min\left\{\frac{\beta_i}{\alpha_i}:\alpha_i \ \mbox{is
even}\right\}}$. Then replace $e$ with $e+1$ and $\beta_j$ with
$\alpha_j+\beta_j$.

Choose continued fraction expansions with odd length $m_i$ and
with $a_2^i,a_4^i,\ldots,a_{m_i-1}^i$ even. The black surface is
orientable in the resulting diagram.

\subsection{The correction term for Seifert fibered spaces}
\label{subsec:seifertlensd}
When $Y$ is the lens space $L(p,q)$ a labelling of $\spinc(Y)$ by
$\zz/p=\{0,1,\ldots p-1\}$ is chosen in \cite[\S 4]{os4}, and the following
recursive formula is given:

$$d(L(p,q),i)=\left(\frac{pq-(2i+1-p-q)^2}{4pq}\right)-d(L(q,r),j),$$
where $i\in\zz/p$ and $r$ and $j$ are the reductions modulo $q$ of
$p$ and $i$ respectively.  (Note Remark \ref{rmk:orient} above
concerning orientation conventions.)

The conjugation action on \spinc structures is given by
$$\bj (i)=q-i-1 \pmod{p},$$
so that the $\bj$-fixed-point-set is
$\displaystyle{\zz\cap\left\{\frac{q-1}{2},\frac{p+q-1}{2}\right\}}$.

\begin{remark}
\label{rmk:froyshov}
It is shown in \cite{s} that the \froyshov\  invariant defined using
Seiberg-Witten theory satisfies the same recursive formula. Therefore a gauge
theoretic version of Theorem \ref{thm:rhomap} based on \cite{f} gives the same
results for lens spaces.
\end{remark}

More generally, if $Y$ is a Seifert fibered rational homology sphere 
$Y(e;(\alpha_1,\beta_1),\ldots,\allowbreak(\alpha_r,\beta_r))$, the following
formula is given in \cite[Corollary 1.5]{os6}:
$$d(Y,\spinct)=\displaystyle{
\mbox{max}\left\{\frac{c_1(\spincs)^2+|G|}{4}:\spincs\in\spinc(X_G),
\spincs|_Y=\spinct\right\}}.$$
Here $G$ is a graph as in Definition \ref{def:seifert} for which
the plumbed manifold $X_G$ is negative definite with $\partial
X_G=Y$, and $|G|$ is the number of vertices of $G$.  This formula
may be interpreted as saying that equality is obtained in
(\ref{thm9.6}) for some $\spincs\in\spinc(X_G)$.  Thus computing
the correction terms for $Y$ is equivalent to computing the $sq$
function on $\cals(X_G)$; in Section \ref{sec:algorithm} we
indicate how to do this for any negative definite four-manifold.

\section{Obstruction algorithm} \label{sec:algorithm}

Given a rational homology sphere $Y$ with the order of $\calh:=H^2(Y;\zz)$ equal to
$h$, and an integer $b\ge 0$ we want to know if $Y$ can bound a negative
definite four-manifold $X$ with $b_2(X)=b$. In view of results of Section
\ref{sec:cohomology} this can be checked in the following sequence of steps:

\begin{itemize}

\item[(1)] consider all factorizations $h=st^2$ with $s,\, t \ge 1$;

\item[(2)] for a fixed factorization, consider all order $t$ subgroups $\calt$
of $\calh$,
and for a fixed $\calt$ consider all order $s$ subgroups $\cals$ of $\calh/\calt$;

\item[(3)] for a fixed $\cals$, consider all negative definite symmetric matrices
$Q$ of rank $b$ that present $\cals$;

\item[(4)] for a fixed $Q$, determine the function $sq:\cals \to \qq$
(see discussion preceding  Theorem \ref{thm:rhomap});

\item[(5)] for all choices of origin in $\calh$ consider all group
monomorphisms $\rho: \cals \to \calh/\calt$, and for a fixed $\rho$
determine the function $d_\rho : \cals \to \qq$;

\item[(6)] if for a particular set of choices above the conclusion of Theorem
\ref{thm:rhomap} holds, then there is no obstruction to $Y$ bounding a negative
definite four-manifold $X$ with $b_2(X)=b$.

\end{itemize}

Note that when $b=0$ the above procedure simplifies significantly (see below
for details). Also, for $b>0$ there is only a finite number of possible choices
in steps (3) and (5); in particular, a complete (but not minimal) set of forms
due to Hermite is described in \cite[Theorem 23]{j}. Similarly, when
determining the function $sq$ in step (4) one can restrict to \spinc structures
whose Chern classes $c$ (modulo torsion) are characteristic vectors in the
hypercube
$$x_i^2 \le c(x_i) < |x_i^2|, \qquad i=1,\ldots,b\, ,$$
where $\{x_i,i=1,\ldots,b\}$ is a basis for $H^2(X,\zz)/\tors$. To
see this note that if the inequality is violated for some $i$,
changing $c$ by an even multiple of the Poincar\'e dual of $x_i$
to make this particular inequality hold, will result in a vector
with no smaller square; moreover, the square only stays the same
if $c(x_i)=|x_i^2|$ (see \cite{os6} for details). A characteristic
vector is the Chern class of a $\bj$-fixed element if and only if
it is in the image of $Q:\zz^b\to\zz^b$. In the rest of this
section we describe in detail the cases $b=0,\ 1$ and $2$ with
emphasis on applications to knots and links.

\subsection{$b=0$} \label{subsec:homball}
A necessary condition for a rational homology sphere $Y$ to bound a rational
homology ball $X$ is that the order of the first homology of $Y$ is a square
(Lemma \ref{lemma:torsion}). The algorithm described above yields the
following.

\begin{proposition}
\label{prop:homball}
Let $X$ be a smooth four-manifold with boundary $Y$, and suppose that
$H_*(X;\qq)\cong H_*(B^4;\qq)$ and the order of $H_1(Y;\zz)$ is $h$; write
$h=t^2$ for some $t\in\nn$. Then there is a spin structure $\spinct_0$ on $Y$
so that
$$d(Y,\spinct_0+\beta)=0 \hbox{\ \ for all\ \ } \beta \in \calt\, ,$$
where $\calt$ denotes the image of $H^2(X;\zz)$ in $H^2(Y;\zz)$.

In particular, if $Y$ is a homology lens space, then given a labelling
$\{\spinct_0,\ldots, \spinct_{h-1}\}$ of \spinc structures on $Y$, there is a
$j_0$ corresponding to a spin structure $\spinct_{j_0}$ so that
$$d(Y,\spinct_{j_0+kt})=0 \hbox{\ \ for all\ \ } k=0,\ldots,t-1.$$
\end{proposition}
\proof
Denote by $\spinct_0$ a \spinc structure on $Y$ that extends to $X$. Then the
set of \spinc structures on $Y$ that extend to $X$ is $\spinct_0+\calt$. Given
that all \spinc structures on the rational homology ball $X$ are torsion,
Theorem \ref{thm:rhomap} implies that $d(Y,\spinct_0+\beta) \ge 0$ for all
$\beta$. Finally, changing the orientation of $X$ and using the fact that the
correction term changes its sign under this operation, gives the other
inequality.
\endpf

Let $K$ be a knot in $S^3$ with branched double cover $Y$.  From the discussion
in Section \ref{sec:knots} we see that if $K$ is slice, then $Y$ bounds a
rational homology ball, and therefore satisfies the conclusion of Proposition
\ref{prop:homball}.

\subsection{$b=1$} \label{subsec:rank1}
When $b_2(X)=1$, the intersection form $Q_X$ of a negative definite manifold
$X$ is represented by $[-s]$, where $h=st^2$ is the order of the first homology
of $Y=\partial X$. Note that in this case $\cals\cong\zz/s$ is cyclic, and
so $Y$ can only bound such an $X$ if $\calh/\calt$ contains a cyclic subgroup of
order $s$ (see discussion preceding Theorem \ref{thm:rhomap} for notation).

Characteristic vectors in $H^2(X;\zz)/\tors$ are given by numbers $x \in \zz$
with the same parity as $s$. A set of \spinc structures on $X$ with maximal
square in their equivalence class in $\cals$ is given by $\spincs_i$,
$i=0,\ldots,s-1$, where the image of $\spincs_i$ modulo torsion is $x_i=2i-s$,
and its square is
$$sq(i)=c_1(\spincs_i)^2=-\frac{(2i-s)^2}{s}\, .$$
Note that $x_0=-s$ corresponds to a $\bj$-fixed element in $\cals$; in case $s$
is even, $x_{s/2}=0$ also gives a $\bj$-fixed element.

Let $L$ be a two component link in $S^3$ with branched double cover $Y$. If the
signature $\sigma(L)=-1$, then according to Murasugi's result $L$ may bound a
cylinder in the four-ball. If this is the case, then $Y$ bounds a negative
definite four-manifold $X$ with $b_2(X)=1$ (see Section \ref{sec:knots}). We
may use the above algorithm to check if this is possible.

\subsection{$b=2$} \label{subsec:rank2}
We now suppose a rational homology sphere $Y$ bounds a negative definite
four-manifold $X$ with $b_2(X)=2$. We denote the order of the first homology of
$Y$ by $h$, and fix a factorization $h=st^2$, where $s$ is the determinant of
the intersection pairing $Q_X$ of $X$. Note that in this case
$\cals\cong\zz^2/Q\zz^2$ has at most two exponents, which puts a homological
restriction on $Y$ bounding such a manifold.

The following classification theorem for rank two quadratic forms is a modified
version of \cite[Theorem 76]{j}.

\begin{theorem}
\label{thm:jones}
Any negative definite form with integer coefficients of rank two and determinant $r>0$ is
equivalent to a reduced form
$$\left[\begin{array}{cc}
a&b\\
b&c
\end{array}\right],$$
with $0\ge2b\ge a\ge c$.  (It follows that $|a|\le2\sqrt{r/3}$.)
\end{theorem}

Let $Q$ be a reduced form as above with coefficients $a,b,c$. Note that this
form presents $\zz/{e_1}\oplus \zz/{e_2}$, where $e_2|e_1$ and
$e_1e_2=r$, if and only if $\gcd(a,b,c)=e_2$. Denote by $\Omega_Q$ the set of
points $(x,y)$ in the plane satisfying the following conditions:
\begin{eqnarray*}
a\le&\!\!\! x\!\!\! &<|a|\, ,\\
c\le &\!\!\! y\!\!\! &<|c|\, ,\\
a-2b+c\le  x&\!\!\!\! - \!\!\!\! & y <|a-2b+c|\, .
\end{eqnarray*}
Call a lattice point $(x,y)\in\zz^2$ {\em characteristic} if $x\equiv a,\
y\equiv c\pmod{2}$.

\begin{proposition}
\label{prop:polygon} Fix a basis for $H_2(X;\zz)/\tors$ so that the matrix
representative $Q$ of the intersection pairing $Q_X$ is reduced. Then the
characteristic lattice points in $\Omega_Q$ are the (images of the) first Chern
classes of a set of \spinc structures on $X$ that have maximal square in their
class in $\cals$. Moreover, any characteristic vector among
$$ (0,0), \quad (a,b), \quad (b,c), \quad (a-b,b-c) $$
gives rise to a $\bj$-fixed element of $\cals$.
\end{proposition}

\proof
Note that two characteristic points correspond to \spinc structures with
isomorphic restrictions to $Y$ if and only if they differ by $2m(a,b)+2n(b,c),$
for some integers $m,n$. A complete set of characteristic representatives is
given by the parallelogram with vertices $\pm(a+b,b+c),\pm(a-b,b-c)$ (taking
all characteristic points in the interior and those in one component of the
boundary minus $\pm(a-b,b-c)$). Observe that each of these points is equivalent
to exactly one characteristic point in $\Omega_Q$. It therefore remains to show
that the corresponding \spinc structures have maximal square in their
equivalence class. Recall that the cup product pairing on $H^2(X;\zz)/\tors$ is
well defined over $\qq$, and its matrix with respect to the Hom-dual basis is
$Q^{-1}$.

As observed after the description of the algorithm, we need only
consider points in the rectangle $\{(x,y)|a\le x < |a|,c\le y <
|c|\}$. Thus it only remains to choose the point with larger
square from any equivalence class having more than one
characteristic point in the rectangle.  This is done by
eliminating points inside the triangles cut out of the rectangle
by the lines $x-y=\pm|a-2b+c|$.
\endpf

It follows from Proposition \ref{prop:polygon} that the numbers
$sq(\alpha)$ (for $\alpha\in\cals$) from Theorem \ref{thm:rhomap}
are given, as an unordered set, by the squares with respect to
$Q^{-1}$ of characteristic points $(x,y)\in\Omega_Q$.  It remains
to order these points with respect to the group structure on
$\cals\cong\zz^2/Q\zz^2\cong\zz/{e_1}\oplus \zz/{e_2}$. The point
$(x_{0,0},y_{0,0})$ may be chosen arbitrarily; for convenience we
choose it to be $\bj$-fixed.  Then choose
$(x_{1,0},y_{1,0})$ so that $\delta_1=
\frac{1}{2}(x_{1,0}-x_{0,0},y_{1,0}-y_{0,0})$ has order $e_1$.
Finally choose $(x_{0,1},y_{0,1})$ so that
$\delta_2= \frac{1}{2}(x_{0,1}-x_{0,0},y_{0,1}-y_{0,0})$ has order
$e_2$ and the subgroups of $\cals$ generated by $\delta_1$ and
$\delta_2$ have trivial intersection. These choices determine the
ordering of the remaining points: $(x_{i,j},y_{i,j})$ is the
unique characteristic point in $\Omega_Q$ with
$$\frac{1}{2}(x_{i,j}-x_{0,0},y_{i,j}-y_{0,0})=i\delta_1+j\delta_2+m(a,b)+n(b,c),\quad m,n\in\zz\, ,$$
and
$$sq(i,j)=(x_{i,j}\,\, y_{i,j})Q^{-1}(x_{i,j}\,\, y_{i,j})^T\, ,$$
for $i=0,\ldots,e_1-1$ and $j=0,\ldots,e_2-1$.

Suppose that $K$ is a knot in $S^3$ with signature $-2$ and branched double
cover $Y$.
From Section \ref{sec:knots} we know that if $g^*(K)=1$, then $Y$ bounds a
four-manifold with $b_2=2$ as above, and we may use the algorithm described at
the beginning of this section to seek a contradiction.  We may similarly get
an obstruction to a three component link with signature $-2$ bounding a genus
zero slice surface.

\section{Examples}
\label{sec:examples}
In this section we list examples of knots and links for which our obstruction
shows that inequality (\ref{eqn:murasugi}) is strict.  We begin with a proof
that the unknotting number of the knot $10_{145}$ is 2.  We list two-bridge
examples in \ref{subsec:2bex} and Montesinos examples in \ref{subsec:montex}.

\subsection{Unknotting number of $10_{145}$}
\label{subsec:10_145}

\begin{figure}[htbp]
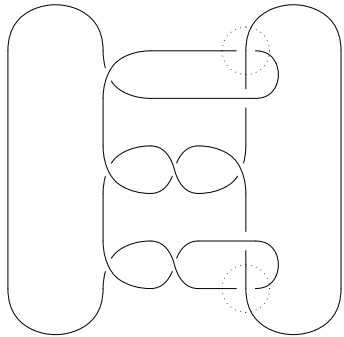

\begin{center}
\ifpic
\leavevmode
\xygraph{
!{0;/r1.5pc/:}
!{\vcap[2]}
!{\xcapv[5]@(0)}[dddd]
!{\vcap[-2]}[u(5)r(5)]
!{\vcap[2]}[rr]
!{\xcapv[5]@(0)}[d(4)l(2)]
!{\vcap[-2]}[u(5)l(3)]
!{\hover}
!{\xcaph[1.8]@(0)}[l(1)d]
!{\xcaph[2.2]@(0)}[ur]
!{\xcapv[0.8]@(0)}[u(1)r(0.2)]
!{\hcap}[l(3.2)d]
!{\xcapv@(0)}[r(3)u(0.8)]
!{\xcapv[0.8]@(0)}[l(3)u(0.2)]
!{\hunder}
!{\htwist}
!{\hunder-}[dlll]
!{\xcapv@(0)}[rrru]
!{\xcapv[0.8]@(0)}[lll]
!{\hunder}
!{\htwist}
!{\xcaph[1.2]@(0)}[ld]
!{\xcaph[0.8]@(0)}[u(0.8)]
!{\xcapv[0.8]@(0)}[u(1.2)r(0.2)]
!{\hcap}[d(1)l(0.2)]
!{\ellipse(0.5){.}}[u(5)]
!{\ellipse(0.5){.}}
}
\else \vskip 5cm \fi
\begin{narrow}{0.4in}{0.2in}
\caption{
{\bf{The Montesinos knot M(1;(3,1),(3,1),(5,2)), or $10_{145}$.}}
Note that changing the circled crossings will give the unknot.}
\label{fig:10_145}
\end{narrow}
\end{center}
\end{figure}

The knot $10_{145}$ in the Rolfsen table is the Montesinos knot
$M(1;(3,1),(3,1),(5,2))$.
From Figure \ref{fig:10_145} we see that the unknotting number is at most 2.
This knot has signature 2 and determinant 3.  Its branched double cover is the
Seifert fibered space $Y(-1;(3,1),(3,1),(5,2))$.  We will show that $-Y$ cannot
bound a negative definite 4-manifold with $b_2=2$.  The correction terms are
$$d(-Y)=\left\{{\displaystyle -\frac{3}{2},-\frac{1}{6},-\frac{1}{6}}\right\}.$$
There are two reduced negative definite forms of rank 2 and determinant 3, namely
the diagonal form
$\left(\begin{array}{cc}
-1 & 0\\
0 & -3 \end{array} \right),$ and
$\left(\begin{array}{cc}
-2 & -1\\
-1 & -2 \end{array} \right)$.
For the first the region $\Omega_Q$ of Proposition \ref{prop:polygon} yields
optimal \spinc structures with squares
$\left\{{\displaystyle -4,-\frac{4}{3},-\frac{4}{3}}\right\}$,
and for the second,
$\left\{{\displaystyle 0,-\frac{8}{3},-\frac{8}{3}}\right\}$.  In
either case there is clearly no map
$$\rho:\zz/3\to\zz/3$$
which satisfies
$$c_1(\spincs_i)^2+2\le4d(-Y,\spinct_{\rho(i)}).$$
It follows that $g^*(10_{145})>1$.  Since the unknotting number is bounded below
by the four-ball genus, we conclude that $g^*=u=2$.

Finally we note that the spanning surface described in \ref{subsec:spanning}
yields the Seifert matrix
$$M=\left(\begin{array}{cccc}
1 & -1 & -1 & 0 \\
0 & 1 & -1 & 0 \\
0 & 0 & 1 & -1 \\
0 & 0 & 0 & 1 \end{array} \right).$$
The vector $x=(1,1,1,0)$ satisfies $xMx^T=0$.  It follows that the Taylor
invariant $m(10_{145})$ (which is the optimal lower bound for $g^*$ from a Seifert
matrix) is 1.

\subsection{A non-cyclic example}
\label{subsec:2cycex}
The Montesinos knot $M(1;(5,2),(5,2),(5,2))$ has signature 4 and
determinant 25.  Its branched double cover $Y=Y(-1;(5,2),(5,2),(5,2))$ has
$H^2(Y)\cong\zz/5\oplus\zz/5$.  The correction terms are

$$d(-Y)=\frac{1}{5}\left\{\begin{array}{ccccc}
 -5 & 1 & -1 & -1 & 1 \\
 -7 & -3 & -7 & 1 & 1 \\
 -3 & -1 & -7 & -1 & -3 \\
 -3 & -3 & -1 & -7 & -1 \\
 -7 & 1 & 1 & -7 & -3
\end{array} \right\},$$
where the array structure indicates the $\zz/5\oplus\zz/5$ action on $\spinc(-Y)$.

Suppose that $-Y$ bounds a negative definite manifold $X$ with $b_2(X)=4$.  Then by
Lemma \ref{lemma:torsion} the intersection pairing of $X$ is either unimodular
or has determinant 25.  If unimodular then it is equivalent to $-I$ (see for
example \cite[Corollary 23]{j}); in this case $\cals$ contains just 1 element
with maximal square $-4$.  The inequality in Theorem \ref{thm:rhomap} now
simply becomes $0\le d_\rho(\alpha)$; however, the correction term of the
$\bj$-fixed element is $-1$.

Now suppose that $Q_X$ has determinant 25.  Note that there are 6 nonnegative
correction terms in the above array.  There are 3 Hermite-reduced negative
definite rank 4 forms
with determinant 25 which present $\zz/5\oplus\zz/5$.  Each of these gives at
least 10 elements $\alpha\in\cals$ with $sq(\alpha)+4\ge0$.  It follows
from Theorem \ref{thm:rhomap} that $-Y$ cannot bound these forms.

This implies $g^*(M(1;(5,2),(5,2),(5,2)))>2$.  From the knot diagram as in
Figure \ref{fig:mont} it is easy to see that the unknotting number is at most
3; thus $g^*=u=3$.

\subsection{Two-bridge examples}
\label{subsec:2bex}
We start with the question of slice two-bridge knots.  Recall a knot $K$ is slice
if $g^*(K)$=0.  It is called {\it ribbon} if it bounds a smoothly immersed
disk in $S^3$ whose singularities come from identifying spanning arcs in $D^2$
with interior arcs in $D^2$.  Ribbon implies slice, however, it is unknown whether
every slice knot is ribbon.

Any slice 2-bridge knot
$S(p,q)$ must have $p=t^2$ from Lemma \ref{lemma:torsion}.  A set of values
of $t$ and $q$ for which $S(t^2,q)$ is ribbon is given in \cite{cg2}.  Using
the Atiyah-Singer $G$-signature theorem, Casson and Gordon \cite{cg2} defined an
invariant which detects when a two-bridge knot is not ribbon
and showed that the known ribbon two-bridge
knots provide the only ribbon examples $S(t^2,q)$ with $t\le 105$.

Fintushel and Stern showed in \cite{fs} that the Casson-Gordon invariant is
equal to an invariant they defined using Yang-Mills theory, and also
showed the invariant detects when a knot is not slice.

The obstruction algorithm described in Subsection \ref{subsec:homball} seems
to give the same results as Casson-Gordon and Fintushel-Stern; we have
verified this for $t\le105$.

Table \ref{table:2bridge}
lists all two-bridge knots and links $S(p,q)$ with $p\le 120$
and $1\le|\sigma|\le4$ for which the obstruction algorithm
shows that inequality
(\ref{eqn:murasugi}) is strict.

\begin{table}[h]
\begin{center}
\begin{tabular}{|l|r|c|c|} \hline
\hskip15pt Link & $\sigma\hskip4pt$ & $m$ & $g^*>$ \\ \hline
$S(60,23)$ & $1$ &  & 0 \\
$S(66,25)$ & $1$ &  & 0 \\
$S(67,39)$ & $2$ & 1 & 1 \\
$S(86,33)$ & $1$ &  & 0 \\
$S(91,53)$ & $2$ & 1 & 1 \\
$S(92,33)$ & $-1$ &  & 0 \\
$S(92,39)$ & $1$ &  & 0 \\
$S(107,28)$ & $-2$ & 1 & 1 \\
$S(107,42)$ & $2$ & 1 & 1 \\
$S(112,43)$ & $1$ &  & 0 \\
$S(114,25)$ & $1$ &  & 0 \\
$S(115,37)$ & $2$ & 1 & 1 \\
$S(115,67)$ & $2$ & 1 & 1 \\
$S(115,87)$ & $-2$ & 1 & 1 \\
\hline
\end{tabular}
\begin{narrow}{0.4in}{0.2in}
\caption{
{\bf{Genus bounds for two-bridge links.}}
Here $\sigma$ is the signature and $m$ is Taylor's invariant.}
\label{table:2bridge}
\end{narrow}
\end{center}
\end{table}

Finally we note that the knots $S(187,101)$ and $S(187,117)$ have the same
Alexander polynomials and Taylor invariants.  The latter has $g^*=1$, but our
algorithm can be used to show that the former has $g^*=2$.

\subsection{More Montesinos examples}
\label{subsec:montex}
Table \ref{table:mont} contains obstructed Montesinos links
$M(e;(\alpha_1,\beta_1),(\alpha_2,\beta_2),(\alpha_3,\beta_3)$ with $-2\le e
\le 1$, $\alpha_i\le 5$, and $|\sigma|\le 4$.  We have also restricted to
links with determinant less than 150.

\vfill\pagebreak

\begin{table}[h]
\begin{center}
\begin{tabular}{|l|c|r|c|c|c|} \hline
\hskip55pt Link & $\mu$ & $\sigma\hskip4pt$ & $H_1(Y)$ & $m$ & $g^*>$ \\ \hline
$M(-2;(3,1),(3,1),(5,3))$ & 1 & $2$ & $\zz/147$ &  & 1 \\
$M(-1;(2,1),(2,1),(5,2))$ & 2 & $-1$ & $\zz/48$ &  & 0 \\
$M(-1;(2,1),(5,2),(5,4))$ & 1 & $2$ & $\zz/135$ &  & 1 \\
$M(-1;(2,1),(5,3),(5,3))$ & 1 & $-2$ & $\zz/135$ &  & 1 \\
$M(-1;(3,1),(3,1),(5,1))$ & 2 & $3$ & $\zz/84$ &  & 1 \\
$M(-1;(3,1),(4,1),(5,4))$ & 1 & $2$ & $\zz/143$ &  & 1 \\
$M(-1;(3,2),(3,2),(5,2))$ & 1 & $2$ & $\zz/123$ &  & 1 \\
$M(-1;(3,2),(4,1),(5,2))$ & 1 & $-2$ & $\zz/139$ &  & 1 \\
$M(0;(2,1),(2,1),(3,2))$ & 2 & $-1$ & $\zz/20$ &  & 0 \\
$M(0;(3,1),(3,1),(5,4))$ & 2 & $1$ & $\zz/66$ &  & 0 \\
$M(0;(3,1),(5,1),(5,4))$ & 2 & $1$ & $\zz/100$ &  & 0 \\
$M(0;(3,1),(5,2),(5,3))$ & 2 & $1$ & $\zz/100$ &  & 0 \\
$M(0;(3,2),(3,2),(5,2))$ & 2 & $1$ & $\zz/78$ &  & 0 \\
$M(0;(3,2),(3,2),(5,4))$ & 2 & $-1$ & $\zz/96$ &  & 0 \\
$M(0;(3,2),(5,1),(5,1))$ & 2 & $-1$ & $\zz/80$ &  & 0 \\
$M(1;(3,1),(3,1),(5,2))$ & 1 & $2$ & $\zz/3$ & 1 & 1 \\
$M(1;(3,1),(5,2),(5,2))$ & 2 & $3$ & $\zz/10$ &  & 1 \\
$M(1;(3,1),(5,4),(5,4))$ & 2 & $-1$ & $\zz/70$ &  & 0 \\
$M(1;(4,1),(4,1),(5,4))$ & 2 & $1$ & $\zz/24$ &  & 0 \\
$M(1;(5,2),(5,2),(5,2))$ & 1 & $4$ & $\zz/5\oplus\zz/5$ &  & 2 \\
\hline
\end{tabular}
\begin{narrow}{0.4in}{0.2in}
\caption{
{\bf{Genus bounds for Montesinos links.}}
Here $\mu$ is the number of components and $Y$ is the branched double cover.}
\label{table:mont}
\end{narrow}
\end{center}
\end{table}

\begin{remark}
\label{rmk:3comp} The reflection of
$M(e;(\alpha_1,\beta_1),\ldots,\allowbreak(\alpha_r,\beta_r))$ is
$M(r-e;(\alpha_1,\alpha_1-\beta_1),\ldots,\allowbreak(\alpha_r,\alpha_r-\beta_r))$.
The four-ball genus of a knot is equal to that of its reflection;
however the same is not true for links, whose signature and four-ball genus
depend on a choice of orientation.  For example the 3-component
link $M(5;(2,1),(2,1),(2,1))$, oriented as in Subsection
\ref{subsec:spanning}, has signature $-2$ and is shown by our
algorithm to have nonzero four-ball genus. Its reflection
$M(-2;(2,1),(2,1),(2,1))$ also has signature $-2$, but the
algorithm yields no information.
\end{remark}


\noindent
Brendan Owens:\\
Dept. of Mathematics and Statistics, McMaster University\\
1280 Main Street West\\
Hamilton, Ontario, Canada L8S 4K1\\
E-mail: owensb@math.mcmaster.ca
\medskip

\noindent
Sa\v{s}o Strle:\\
Dept. of Mathematics and Statistics, McMaster University\\
1280 Main Street West\\
Hamilton, Ontario, Canada L8S 4K1\\
E-mail: strles@math.mcmaster.ca


\begin{thebibliography}{99}
\bibitem{b} G.~Burde. \textsl{Uber das Geschlecht und die Faserbarkeit von
    Montesinos-knoten}, Abh. Math. Sem. Univ. Hamburg {\bf 54} 1984, 199--226.
\bibitem{bz} G.~Burde \& H.~Zieschang. \textsl{Knots},
    Walter de Gruyter, 1985.
\bibitem{cg2} A.~J.~Casson, C.~McA.~Gordon. \textsl{Cobordism of classical knots},
    \`A la recherche de la topologie perdue, Progr. Math. {\bf 62} 1986, 181--199.
\bibitem{cg} A.~J.~Casson, C.~McA.~Gordon. \textsl{On slice knots in dimension
    three},
    Algebraic and geometric topology, Part 2, pp. 39--53,
    Proc. Sympos. Pure Math., XXXII, Amer. Math. Soc., Providence, R.I., 1978.
\bibitem{fs} R.~Fintushel \& R.~J.~Stern. \textsl{Rational homology cobordisms
    of spherical space forms}, Topology {\bf 26} 1987, 385--393.
\bibitem{f} K.~A.~Fr{\o}yshov. \textsl{The Seiberg-Witten equations and
    four-manifolds with boundary}, Math.~Res.~Lett. {\bf 3}
    1996, 373--390.
\bibitem{hnk} F.~Hirzebruch, W.~D.~Neumann \& S.~S.~Koh. \textsl{Differentiable
    manifolds and quadratic forms},
    Lecture Notes in Pure and Applied Math. {\bf 4}, Marcel Dekker, 1971.
\bibitem{j} B.~W.~Jones. \textsl{The arithmetic theory of quadratic forms},
    Math.~Assoc.~of America, 1950.
\bibitem{kt} L.~H.~Kauffman \& L.~R.~Taylor. \textsl{Signature of links},
Trans.~Amer.~Math.~Soc. {\bf 216} 1976, 351--365.
\bibitem{m} K.~Murasugi. \textsl{On a certain numerical invariant of link types},
Trans.~Amer. Math.~Soc. {\bf 117}  1965, 387--422.
\bibitem{mont} J.~M.~Montesinos. \textsl{Seifert manifolds that are ramified
two-sheeted cyclic coverings},
Bol.~Soc.~Mat.~Mexicana (2) {\bf 18} 1973, 1--32.
\bibitem{o} P.~Orlik. \textsl{Seifert manifolds},
    Lecture Notes in Math. {\bf 291}, 1972.
\bibitem{os4} P.~\ozsvath \& Z.~Szab\'{o}. \textsl{Absolutely graded Floer
    homologies and intersection forms for four-manifolds with boundary},
    math.SG/0110170, 2002.
\bibitem{os6} P.~\ozsvath \& Z.~Szab\'{o}. \textsl{On the Floer
    homology of plumbed three-manifolds},
    math.SG/0203265, 2002.
\bibitem{os11} P.~\ozsvath \& Z.~Szab\'{o}. \textsl{Knot Floer
    homology and the four-ball genus},
    math.GT/0301149, 2003.
\bibitem{s} S.~Strle. \textsl{In preparation}
\bibitem{t} L.~R.~Taylor. \textsl{On the genera of knots},
    Lecture Notes in Math. {\bf 722} 1972, 144--154.
\end{thebibliography}
\end{document}